\documentclass[a4paper]{article}

\usepackage[inner=30mm,outer=30mm,textheight=225mm]{geometry}
\usepackage{amsmath}
\usepackage{amssymb}
\usepackage{array}
\usepackage{enumerate}
\usepackage{footnote}
\usepackage{graphicx}
\usepackage{psfrag}
\usepackage{theorem}

\newtheorem{theorem}{Theorem}[section]

\newtheorem{defn}[theorem]{Definition}

\newtheorem{conjecture}[theorem]{Conjecture}

\newcommand{\qed}{\hfill\rule[-0.5mm]{1.5mm}{3.0mm}}
\newtheorem{proofthm}{Proof}

\newcommand{\annref}{\bar{A}}

\newcommand{\discref}{\bar{D}}
\renewcommand{\epsilon}{\varepsilon}

\newcommand{\kb}{K^2}

\newcommand{\mb}{M^2}
\newcommand{\mobius}{M\"{o}bius}

\newcommand{\poly}[1]{$\mathit{#1}$}
\newcommand{\polyedge}[1]{$\overline{\mathit{#1}}$}
\newcommand{\polysq}[1]{$\square \mathit{#1}$}
\newcommand{\polytri}[1]{$\Delta \mathit{#1}$}
\newcommand{\ppirr}{{$\mathbb{P}^2$-ir\-re\-du\-ci\-ble}}

\newcommand{\regina}{{\em Regina}}

\newcommand{\rpp}{\mathbb{R}P^2}
\newcommand{\rps}{\mathbb{R}P^3}

\newcommand{\sfs}[2]{\mathrm{SFS}\left(#1: #2\right)}
\newcommand{\sfslong}{Seifert fibred space}

\newcommand{\snappea}{{\em SnapPea}}

\newcommand{\sss}{S^3}

\newcommand{\torus}{T^2}
\newcommand{\twisted}{\stackrel{\smash{\protect\raisebox{-1mm}[0pt][0pt]{$\scriptstyle\sim$}}}{\times}}
\newcommand{\twt}{\tilde{T}}

\newcommand{\Z}{\mathbb{Z}}

\newcommand{\fcaption}[1]{\caption{#1}}
\newcommand{\tcaption}[1]{\caption{#1}}

\newcommand{\homtwo}[4]{
    \mbox{\scriptsize \renewcommand{\arraystretch}{1}
        $\! \left[ \begin{array}{@{\ }c@{\ }c@{\ }} #1 & #2 \\ #3 & #4
        \end{array} \right]$
    }
}

\newcommand{\twolines}[3]{
    \begin{tabular}{#1} #2 \\ #3 \end{tabular}
}

\title{Observations from the 8-tetrahedron \\ non-orientable census}
\author{Benjamin A.~Burton}

\date{September 15, 2005}

\begin{document}

\maketitle

\abstract{Through computer enumeration with the aid of topological
results, we catalogue all 18 closed non-orientable {\ppirr} 3-manifolds
that can be formed from at most eight tetrahedra.  In addition we give
an overview as to how the 100 resulting minimal triangulations are
constructed.  Observations and conjectures are drawn from the census
data, and future potential for the non-orientable census is discussed.
Some preliminary nine-tetrahedron results are also included.}

\section{Introduction} \label{s-intro}

Several surveys have been performed in recent years of all small 3-manifold
triangulations satisfying particular properties.  One of the key strengths
of such a census is in examining {\em minimal triangulations}
(triangulations of 3-manifolds that use as few tetrahedra as possible).

Minimal triangulations are still poorly understood.  Many necessary
conditions for minimality can be found in the literature (see
\cite{burton-facegraphs,0-efficiency,italian9,matveev6} for some
examples).  However, sufficient conditions are much more difficult to
find.  Most of the positive results regarding minimality rely on
exhaustive censuses such as these.

Beyond its use in studying minimal triangulations, a census
also forms a useful body of examples for testing
conjectures and searching for patterns.  Section~\ref{s-obs} illustrates
some conjectures arising from the non-orientable census described in
this paper.

We restrict our attention here to closed {\ppirr} 3-manifolds.
Examples of other censuses involving manifolds with boundaries or cusps can
be seen in the results of Callahan, Hildebrand and Weeks \cite{cuspedcensus}
and Frigerio, Martelli and Petronio \cite{italian-heegaard}.

The extent of census data known to date for minimal 3-manifold
triangulations is fairly small.  This is due to the
computational difficulty of performing such a survey --- the number of
potential triangulations to examine grows worse than exponentially
with the number of tetrahedra.

Closed orientable 3-manifolds have been
surveyed successively by Matveev for six tetrahedra \cite{matveev6},
Ovchinnikov for seven tetrahedra, Martelli and Petronio for nine tetrahedra
\cite{italian9}, and recently Martelli for ten tetrahedra \cite{italian10}.

Closed non-orientable 3-manifolds are less well studied.  The
six-tetrahedron and seven-tetra\-hedron cases were tackled independently
by Amendola and Martelli \cite{italian-nor6,italian-nor7} and by Burton
\cite{burton-nor7}.  Only eight different non-orientable
3-manifolds are found up to
seven tetrahedra (and none at all are found below six tetrahedra).  In this
sense, the seven-tetrahedron results are but a taste of what lies
ahead.

The methods of these different authors are notably distinct.
Amendola and Martelli do not use a direct computer search, but instead
employ more creative techniques.  For seven tetrahedra \cite{italian-nor7}
they examine orientable
double covers and invoke the results of the nine-tetrahedron orientable
census \cite{italian9}.  More remarkable is their six-tetrahedron
census \cite{italian-nor6}, which is purely theoretical and makes no use
of computers at all.

On the other hand, Burton focuses on minimal triangulations
of these eight different 3-manifolds and their combinatorial structures.
With three exceptions in the smallest case (six tetrahedra), the
compositions of the 41 different minimal triangulations found in the
census are described in
detail and generalised into infinite families \cite{burton-nor7}.

The work presented here extends the non-orientable census results to eight
tetrahedra.  Both 3-manifolds and all their minimal triangulations are
enumerated and placed in the context of earlier results.  In total there
are 10 new 3-manifolds with 59 different triangulations.  All 59 of
these minimal triangulations fit within the families described in
\cite{burton-nor7}.

It is worth noting that the list of non-orientable 3-manifolds formed
from $\leq 8$ tetrahedra is equivalently a list of non-orientable
3-manifolds with Matveev complexity $\leq 8$.  Matveev defines the
complexity of a 3-manifold in terms of special spines
\cite{matveev-complexity}, and it is proven by Martelli and Petronio
\cite{italian-decomp} that for all closed {\ppirr} 3-manifolds other
than $\sss$, $\rps$ and $L_{3,1}$, this is
equivalent to the number of tetrahedra in a minimal triangulation.

All of the computational work was performed using {\regina},
a software package that performs a variety of different calculations and
procedures in 3-manifold topology \cite{regina,burton-regina}.
The program \regina, its source code
and accompanying documentation are freely available from
{\tt http://\allowbreak regina.\allowbreak sourceforge.\allowbreak net/}.

In the remainder of Section~\ref{s-intro} we describe in detail the
census parameters and give a concise summary of the results.
Section~\ref{s-tri} presents an overview of how the different minimal
triangulations are constructed, though the reader is referred to
\cite{burton-nor7} for finer details (an appendix is provided to
match the individual census triangulations to the detailed constructions
of \cite{burton-nor7}).  Finally, Section~\ref{s-obs} contains some
observations and conjectures drawn from the census results, and closes
with some remarks regarding future directions of the non-orientable
census.  Partial results from the nine-tetrahedron census (which is
currently under construction) are briefly discussed.

Special thanks must go to J.~Hyam Rubinstein for many helpful
discussions throughout the course of this research.  Thanks are also due
to the University of Melbourne and the Victorian Partnership for
Advanced Computing, both of which have provided computational support
for this and related research.

\subsection{Summary of Results} \label{s-intro-results}

As with the previous closed censuses described above,
we consider only triangulations satisfying the following constraints.
\begin{itemize}
    \item {\em Closed:}  The triangulation is of a closed
    3-manifold.  In particular it has no boundary faces, and each vertex
    link is a 2-sphere.
    \item {\em \ppirr:}  The underlying 3-manifold has no
    embedded two-sided projective planes, and furthermore every
    embedded 2-sphere bounds a ball.
    \item {\em Minimal:}  The underlying 3-manifold cannot be
    triangulated using strictly fewer tetrahedra.
\end{itemize}

Requiring triangulations to be {\ppirr} and minimal keeps the number of
triangulations down to manageable levels, focussing only upon the simplest
triangulations of the simplest 3-manifolds (from which more complex
3-manifolds can be constructed).

The main result of this paper is the following.  As with most censuses
described in the literature, its proof relies upon an exhaustive
computer search.  This search was performed using the software package
\regina, with the help of several results described in \cite{burton-facegraphs}
to increase the efficiency of the search algorithm.  For more details on
how the search algorithm is structured, see the seven-tetrahedron census
paper \cite{burton-nor7}.

\begin{theorem}[Census Results] \label{t-census}
Consider all closed non-orientable {\ppirr} 3-mani\-folds that can be
triangulated using at most eight tetrahedra.  This set contains 18
different 3-manifolds with a total of 100 minimal triangulations between
them, as summarised in Tables~\ref{tab-summary} and~\ref{tab-3mfds}.
\end{theorem}

It should be noted that, when restricted to $\leq 7$ tetrahedra, the
eight different 3-manifolds obtained match precisely with the lists
obtained by Amendola and Martelli \cite{italian-nor7}.

For complete details of the 100 minimal triangulations, a
data file may be downloaded from the {\regina} website
\cite{regina}.\footnote{The eight tetrahedron non-orientable census data
is also bundled with {\regina} version 4.2.1 or later.
It can be found in the {\em File $\rightarrow$ Open Example} menu.}
When opened within {\regina}, the triangulations may be examined in
detail along with various properties of interest such as algebraic
invariants and normal surfaces.

As promised in Theorem~\ref{t-census}, a brief summary of results
appears in Table~\ref{tab-summary}.  Here we see overall totals,
split according to the number of tetrahedra in the minimal
triangulations for each 3-manifold.  Note that
each triangulation is counted once up to isomorphism, i.e., a
relabelling of the tetrahedra within the triangulation
and their individual faces.

\begin{table}[htb]
\tcaption{Summary of closed non-orientable census results}
\begin{center} \begin{tabular}{|r|r|r|}
    \hline
    \bf Tetrahedra & \bf 3-Manifolds & \bf Triangulations \\
    \hline
    $\leq 5$ & 0 & 0 \\
    6 & 5 & 24 \\
    7 & 3 & 17 \\
    8 & 10 & 59 \\
    \hline
    Total & 18 & 100 \\
    \hline
\end{tabular} \end{center}
\label{tab-summary}
\end{table}

Two striking observations can be made from Table~\ref{tab-summary},
which have been made before \cite{italian-nor7,burton-nor7} but are
worth repeating here.
These are that (i) there are no closed non-orientable
{\ppirr} triangulations at all with $\leq 5$ tetrahedra, and that (ii)
the number of minimal triangulations is much larger than the number of
3-manifolds.  Indeed, most 3-manifolds in the census can be realised by
several different minimal triangulations, as seen again in the next
table.

Table~\ref{tab-3mfds} provides finer detail for each of the 18 different
3-manifolds, including the number of minimal triangulations for each
3-manifold and the first homology group.  The notation used for
describing 3-manifolds is as follows.
\begin{itemize}
    \item $\torus \times I / \homtwo{p}{q}{r}{s}$ represents the torus
    bundle over the circle with monodromy $\homtwo{p}{q}{r}{s}$;
    \item $\sfs{B}{\ldots}$ represents a non-orientable {\sfslong} over
    the base orbifold $B$, where $\rpp$ and $\discref$ represent the
    projective plane and the disc with reflector boundary respectively.
    The remaining arguments ($\ldots$) describe the exceptional fibres.
\end{itemize}

\renewcommand{\arraystretch}{1.5}
\begin{table}[htb]
\tcaption{Details for each closed non-orientable {\ppirr} 3-manifold}
\small
\begin{center} \begin{tabular}{|c|l|c|l|l|}
\hline
\bf Tetrahedra & \bf 3-Manifold & \bf Triangulations & \bf Homology \\
\hline
6 & $\torus \times I / \homtwo{1}{1}{1}{0}$ & 1 &
    $\Z$ \\
  & $\torus \times I / \homtwo{0}{1}{1}{0}$ & 6 &
    $\Z \oplus \Z$ \\
  & $\torus \times I / \homtwo{1}{0}{0}{-1}$ & 3 &
    $\Z \oplus \Z \oplus \Z_2$ \\
  & $\sfs{\rpp}{(2,1)\ (2,1)}$ & 9 &
    $\Z \oplus \Z_4$ \\
  & $\sfs{\discref}{(2,1)\ (2,1)}$ & 5 &
    $\Z \oplus \Z_2 \oplus \Z_2$ \\
\hline
7 & $\torus \times I / \homtwo{2}{1}{1}{0}$ & 4 &
    $\Z \oplus \Z_2$ \\
  & $\sfs{\rpp}{(2,1)\ (3,1)}$ & 10 &
    $\Z$ \\
  & $\sfs{\discref}{(2,1)\ (3,1)}$ & 3 &
    $\Z \oplus \Z_2$ \\
\hline
8 & $\torus \times I / \homtwo{3}{1}{1}{0}$ & 10 &
    $\Z \oplus \Z_3$ \\
  & $\torus \times I / \homtwo{3}{2}{2}{1}$ & 2 &
    $\Z \oplus \Z_2 \oplus \Z_2$ \\
  & $\sfs{\rpp}{(2,1)\ (4,1)}$ & 10 &
    $\Z \oplus \Z_2$ \\
  & $\sfs{\rpp}{(2,1)\ (5,2)}$ & 10 &
    $\Z$ \\
  & $\sfs{\rpp}{(3,1)\ (3,1)}$ & 7 &
    $\Z \oplus \Z_6$ \\
  & $\sfs{\rpp}{(3,1)\ (3,2)}$ & 9 &
    $\Z \oplus \Z_3$ \\
  & $\sfs{\discref}{(2,1)\ (4,1)}$ & 3 &
    $\Z \oplus \Z_2 \oplus \Z_2$ \\
  & $\sfs{\discref}{(2,1)\ (5,2)}$ & 3 &
    $\Z \oplus \Z_2$ \\
  & $\sfs{\discref}{(3,1)\ (3,1)}$ & 3 &
    $\Z \oplus \Z_3$ \\
  & $\sfs{\discref}{(3,1)\ (3,2)}$ & 2 &
    $\Z \oplus \Z_3$ \\
\hline
\end{tabular} \end{center}
\label{tab-3mfds}
\end{table}
\renewcommand{\arraystretch}{1}

The most immediate observation is that the eight-tetrahedron census
offers little more variety than the six- and seven-tetrahedron censuses
that came before it.  The census is populated entirely by torus bundles
and by {\sfslong}s over $\rpp$ or $\discref$ with two exceptional fibres.
Preliminary results do suggest that the nine-tetrahedron census will reveal
more variety than this; see Section~\ref{s-obs} for further discussion.

Finally it is worth noting that, as observed by Amendola and Martelli
\cite{italian-nor6,italian-nor7}, all four flat Klein bottle bundles can
be triangulated with only six tetrahedra.  These include all six-tetrahedron
manifolds in the table except for $\torus \times I / \homtwo{1}{1}{1}{0}$.

\section{Constructing Minimal Triangulations} \label{s-tri}

In the seven-tetrahedron census paper \cite{burton-nor7},
the combinatorial structures of the 41 census triangulations
are described in full detail.
A number of parameterised families are presented,
precise parameterised constructions are given for triangulations in these
families, and the resulting 3-manifolds are identified.

Having extended the census to eight tetrahedra,
all of the additional 59 triangulations
are found to belong to these same parameterised families.  We
therefore refer the reader to \cite{burton-nor7} for details of their
construction.  Here we present a simple overview of each family,
showing how their triangulations are pieced
together to form 3-manifolds of various types.  We do go into a little
detail, since these families feature in some of the conjectures
of Section~\ref{s-obs}.

For completeness, the appendix contains a full listing with the
precise parameters for each census triangulation.  This allows the
triangulations to be fully reconstructed and cross-referenced
against \cite{burton-nor7}, though of course the reader is
invited to download the 100 triangulations instead
as a {\regina} data file as described in the introduction.

There are three broad families of triangulations to describe.  These are
the layered surface bundles, the plugged thin $I$-bundles and the
plugged thick $I$-bundles.  Each is discussed in its own section below.

%
%

\subsection{Layered Surface Bundles}

A {\em layered surface bundle} produces either a torus bundle or a Klein
bottle bundle over the circle.  We postpone a formal definition for the
moment, instead giving a broad overview of the construction.

\begin{figure}[htb]
\psfrag{Bdry}{{\small Boundary $F$ (2 faces)}}
\psfrag{Prod}{{\small Product $F \times I$}}
\psfrag{Bdry2}{{\small New boundary $F$ (2 faces)}}
\psfrag{Layer}{{\small Layering to change boundary curves}}
\psfrag{ID}{{\small \twolines{c}{Identification of}{boundaries}}}
\centerline{\includegraphics[scale=0.6]{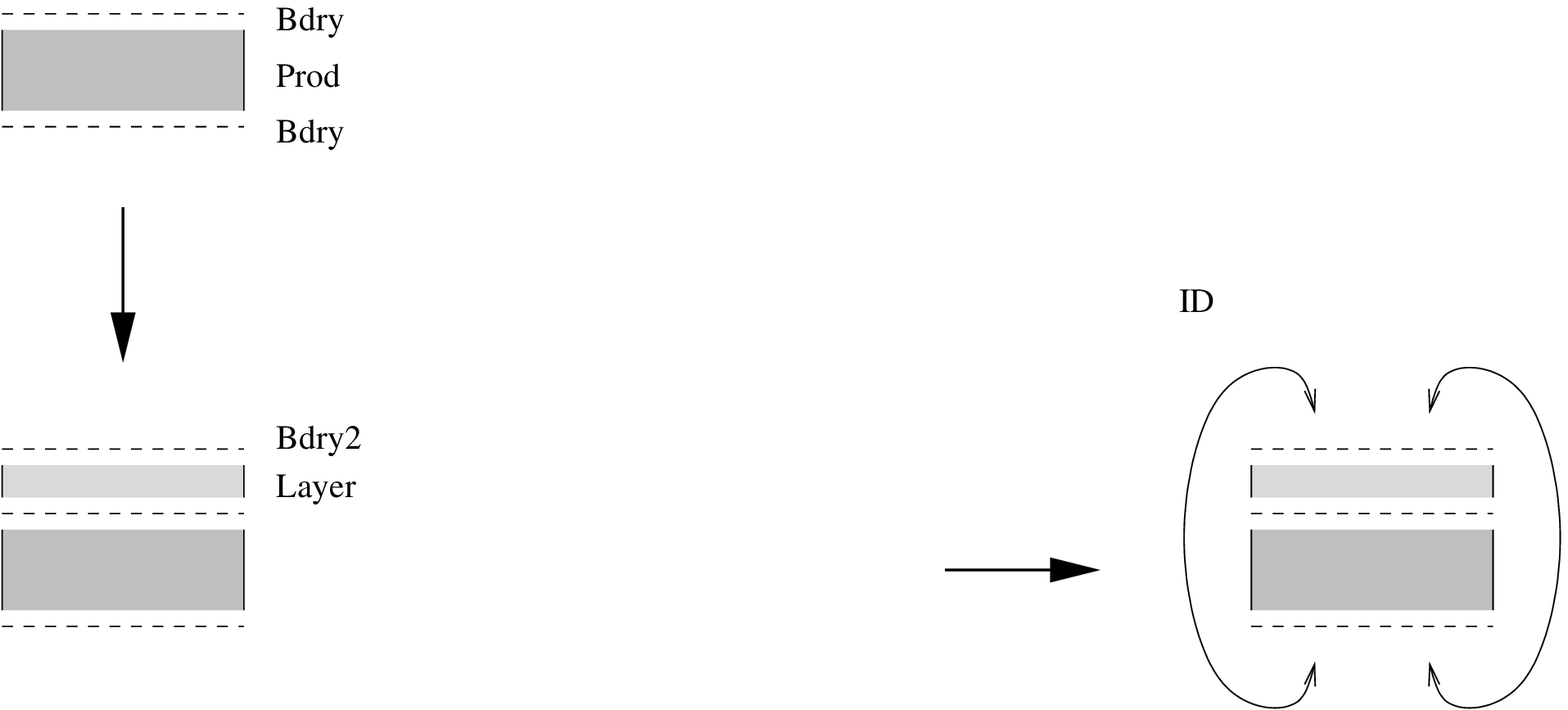}}
\fcaption{Constructing a layered surface bundle}
\label{fig-lsb}
\end{figure}

Figure~\ref{fig-lsb} illustrates the general structure of a layered
surface bundle.  First the product $\torus \times I$ or $\kb \times I$
is constructed (where the surface $F$ in the diagram is either
$\torus$ or $\kb$ for the torus or Klein bottle
accordingly).  This leaves two boundary surfaces,
each of which are then identified according to some specified monodromy.
If this is impossible because the boundary edges do not match, some
additional tetrahedra may be layered onto one of the boundary surfaces
to adjust the boundary edges accordingly.

\subsubsection{Components}

We continue with enough detail to allow a precise definition of a
layered surface bundle as seen in Definition~\ref{d-lsb} below.  This
requires us to describe more precisely how the product $\torus \times I$
or $\kb \times I$ is formed, as well as what a layering entails.

\begin{defn}[Untwisted Thin $I$-Bundle] \label{d-unthinibundle}
    An {\em untwisted thin $I$-bundle} over some closed surface $F$
    is a triangulation of the product $F \times I$ formed as follows.

    Consider the interval $I = [0,1]$.  The product $F \times I$
    is naturally foliated by surfaces $F \times \{x\}$ for
    $x \in [0,1]$.  When restricted to an individual tetrahedron,
    we require that this
    foliation decomposes the tetrahedron into either triangles or
    quadrilaterals as illustrated in Figure~\ref{fig-ibundletet}.
    Note that every vertex lies on one of the boundaries $F \times \{0\}$
    or $F \times \{1\}$, as does the upper face in the triangular case
    and the upper and lower edges in the quadrilateral case.

    In particular, the surface $F \times \{\frac12\}$ meets every
    tetrahedron in precisely one triangle or quadrilateral.  We refer to
    $F \times \{\frac12\}$ as the {\em central surface} of the $I$-bundle.

    \begin{figure}[htb]
    \centerline{\includegraphics[scale=0.6]{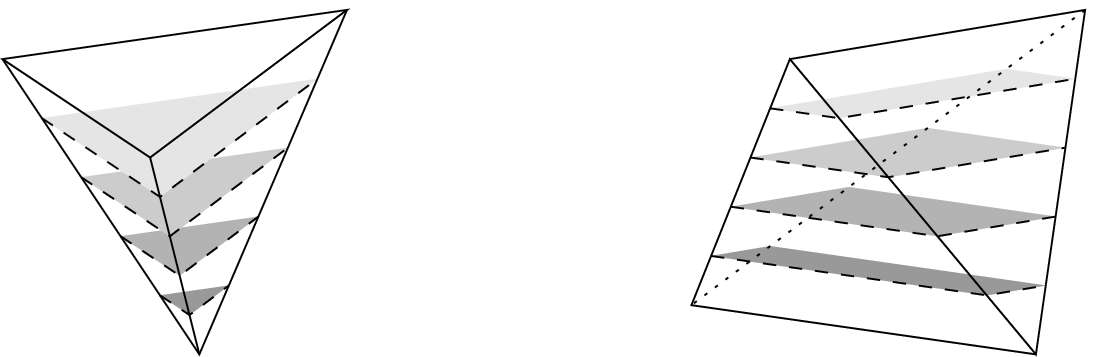}}
    \fcaption{Decomposing a tetrahedron into triangles or quadrilaterals}
    \label{fig-ibundletet}
    \end{figure}
\end{defn}

\begin{figure}[htb]
\psfrag{Prod}{{\small \twolines{c}{The product}{$\torus \times I$}}}
\psfrag{Bdry}{{\small \twolines{c}{The two torus}{boundaries}}}
\psfrag{Central}{{\small The central torus}}
\centerline{\includegraphics[scale=0.6]{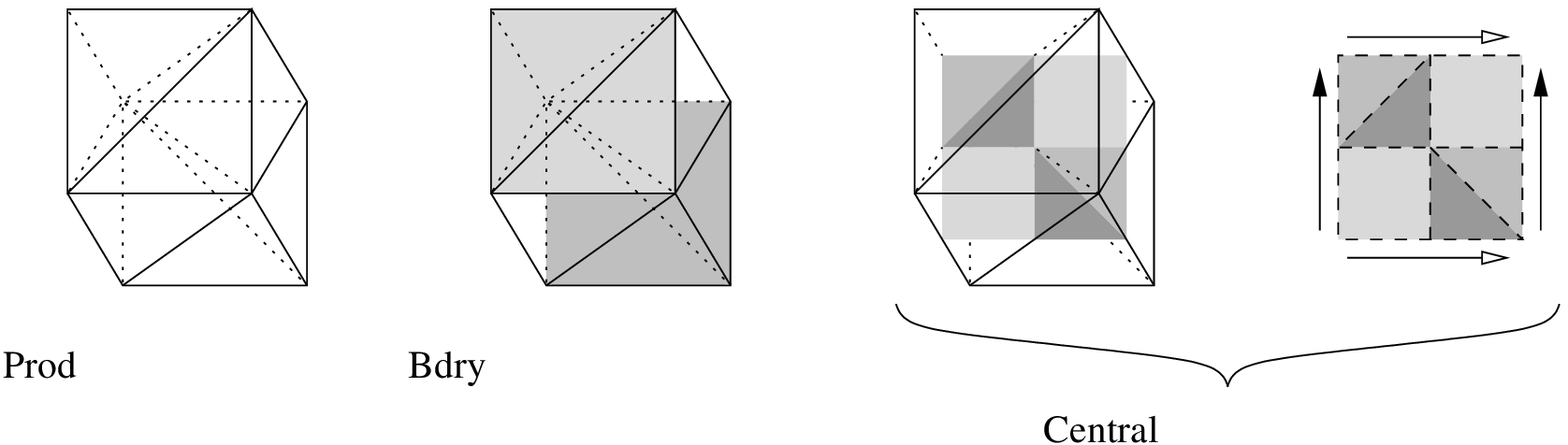}}
\fcaption{An example of a thin $I$-bundle over the torus}
\label{fig-lsb-t62}
\end{figure}

An example of an untwisted thin $I$-bundle over the torus is
illustrated in Figure~\ref{fig-lsb-t62}.  This triangulation consists of
six tetrahedra arranged into a cube.  The front and back faces of the
cube form the boundary tori, which are shaded in the second diagram of the
sequence.  The remaining faces are identified in the usual way for a
torus (the top identified with the bottom and the left identified with
the right).

The central surface $\torus \times \{\frac12\}$ is shown in the
third diagram.  In the fourth diagram we can see precisely how the six
tetrahedra divide this central torus into six cells, each a triangle or
quadrilateral, with the arrows indicating which edges are identified
with which.

It follows from Definition~\ref{d-unthinibundle} that a thin
$I$-bundle is ``only one tetrahedron thick''.  That is, each tetrahedron
runs all the way from one boundary surface to the other, as does each
non-boundary edge.

As a final note, it should be observed that the decomposition of the
central surface offers enough information to completely reconstruct the
thin $I$-bundle.  This is because each triangle or quadrilateral of the
central surface corresponds to one tetrahedron, and the adjacencies of
the triangles and quadrilaterals dictate the corresponding adjacencies
between tetrahedra.

We move now to describe a layering, a well-known procedure
by which a single tetrahedron is attached to a boundary surface in order
to rearrange the boundary edges.

\begin{defn}[Layering] \label{d-layering}
    Consider a triangulation with some boundary component $B$.
    A {\em layering} involves attaching a single tetrahedron $\Delta$
    to the boundary $B$ as follows.  Two adjacent faces of $\Delta$ are
    identified directly with two adjacent faces of $B$, and the
    remaining two faces of $\Delta$ become new boundary faces.  This
    procedure is illustrated in Figure~\ref{fig-layering}.

    \begin{figure}
    \psfrag{e}{{\small $e$}}
    \psfrag{f}{{\small $f$}}
    \psfrag{T}{{\small $\Delta$}}
    \psfrag{B}{{\small $B$}}
    \centerline{\includegraphics[scale=0.7]{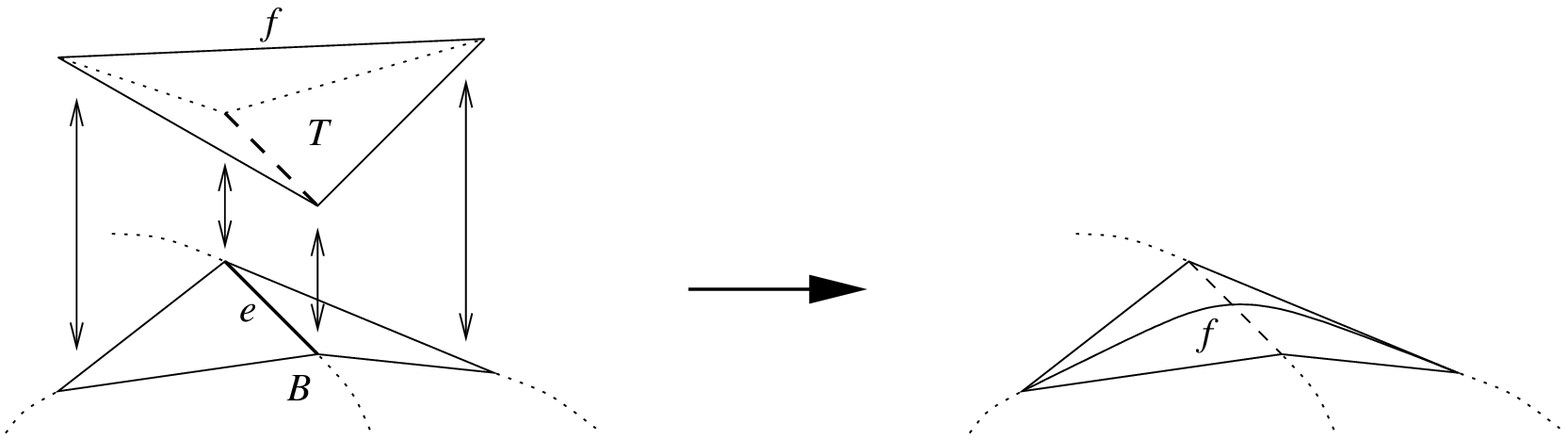}}
    \caption{Performing a layering}
    \label{fig-layering}
    \end{figure}
\end{defn}

The underlying 3-manifold is unchanged --- the primary effect of the
layering is to alter the curves formed by the edges on the boundary.
This is illustrated in the right-hand diagram of
Figure~\ref{fig-layering}, where the old boundary edge $e$ has been made
internal and a new different boundary edge $f$ has appeared in its place.

Given Definitions~\ref{d-unthinibundle} and~\ref{d-layering}, we can
now define a layered surface bundle precisely.

\begin{defn}[Layered Surface Bundle] \label{d-lsb}
    A {\em layered torus bundle} or a {\em layered Klein bottle bundle}
    is a triangulation formed as follows.
    Let $F$ be either the torus or the Klein bottle respectively.
    An untwisted $I$-bundle over $F$ is formed, such that each boundary
    $F \times \{0\}$ and $F \times \{1\}$ consists of precisely two
    faces.  A series of zero or more tetrahedra are then layered onto the
    boundary $F \times \{1\}$, resulting in a new boundary surface $F'$,
    again with precisely two faces.  Finally the surfaces
    $F \times \{0\}$ and $F'$ are identified according to some
    homeomorphism of the original surface $F$.

    For convenience, we refer to both layered torus bundles and layered
    Klein bottle bundles as {\em layered surface bundles}.
\end{defn}

It is clear that the 3-manifold formed from 
a layered torus bundle or a layered Klein bottle bundle
is a torus bundle or Klein bottle bundle over the circle respectively.
Once more the reader is referred to Figure~\ref{fig-lsb} for a pictorial
representation of this procedure.

\subsubsection{Census Triangulations}

The layered surface bundles that appear in the census give rise to the
following 3-manifolds.  From layered torus bundles
we obtain the six manifolds
\renewcommand{\arraystretch}{1.5}
\[ \begin{array}{lll}
\torus \times I / \homtwo{1}{1}{1}{0}, \quad &
\torus \times I / \homtwo{0}{1}{1}{0}, \quad &
\torus \times I / \homtwo{1}{0}{0}{-1}, \\
\torus \times I / \homtwo{2}{1}{1}{0}, \quad &
\torus \times I / \homtwo{3}{1}{1}{0}, \quad &
\torus \times I / \homtwo{3}{2}{2}{1}.
\end{array} \]
From layered Klein bottle bundles we obtain the four flat manifolds
\[ \begin{array}{ll}
\torus \times I / \homtwo{0}{1}{1}{0}, \quad &
\torus \times I / \homtwo{1}{0}{0}{-1}, \\
\sfs{\rpp}{(2,1)\ (2,1)}, \quad &
\sfs{\discref}{(2,1)\ (2,1)},
\end{array} \]
\renewcommand{\arraystretch}{1}
each of which has an alternate expression as a Klein bottle bundle over
the circle.

Table~\ref{tab-lsb} places these observations within the context of the
overall census.  Specifically, it
lists the number of different layered surface bundles that
appear in the census for each number of tetrahedra, as well number of
different 3-manifolds that they describe.  Note that there are only
eight distinct 3-manifolds in total, since in the lists above the torus bundles
$\torus \times I / \homtwo{0}{1}{1}{0}$ and
$\torus \times I / \homtwo{1}{0}{0}{-1}$ each appear twice.

\begin{table}[htb]
\tcaption{Frequencies of layered surface bundles within the census}
\begin{center} \begin{tabular}{|r|rl|rl|}
    \hline
    \bf Tetrahedra & \multicolumn{2}{c|}{\bf 3-Manifolds} &
    \multicolumn{2}{c|}{\bf Triangulations} \\
    \hline
    6 & 5 & (out of 5) & 15 & (out of 24) \\
    7 & 1 & (out of 3) & 4 & (out of 17) \\
    8 & 2 & (out of 10) & 12 & (out of 59) \\
    \hline
\end{tabular} \end{center}
\label{tab-lsb}
\end{table}

Again it can be observed that there are significantly more
triangulations than 3-manifolds.  This is because there are several
different choices for the initial thin $I$-bundle, as well as several
different boundary homeomorphisms (and thus several different layerings)
that can be used to describe the same 3-manifold.

%
%

\subsection{Plugged Thin $I$-Bundles}

\begin{figure}[htb]
\psfrag{TP}{{\small Twisted product $\torus \twisted I$ (6 tets)}}
\psfrag{Bdry}{{\small Boundary torus (4 faces)}}
\psfrag{LST}{{\small \twolines{ll}{Two layered}{solid tori}}}
\centerline{\includegraphics[scale=0.6]{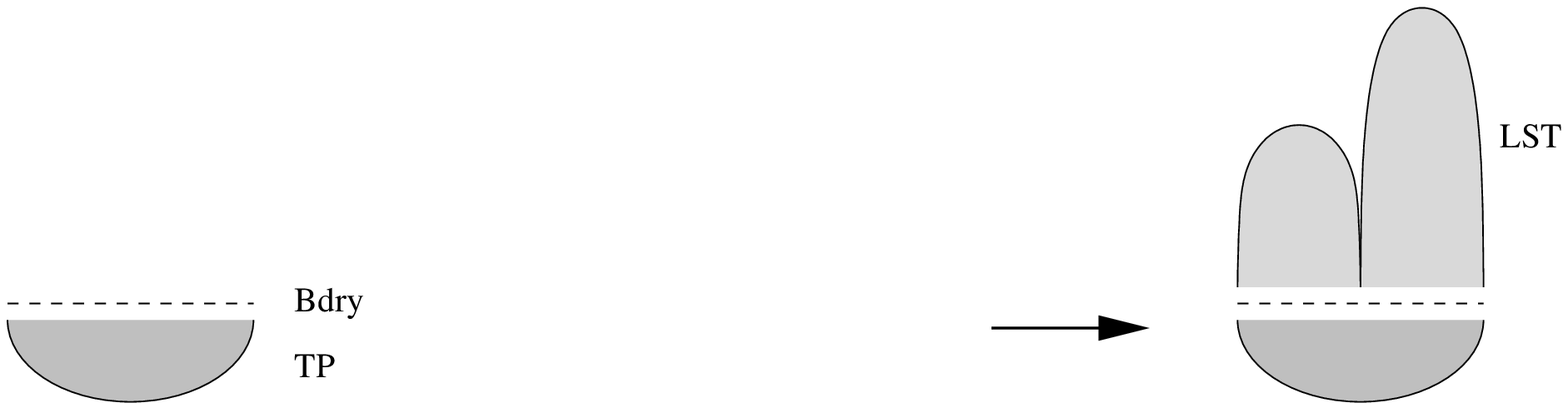}}
\fcaption{Constructing a plugged thin $I$-bundle}
\label{fig-pluggedthin}
\end{figure}

A {\em plugged thin $I$-bundle} allows us to create a non-orientable
{\sfslong} with two exceptional fibres.  It begins with a
six-tetrahedron triangulation of the twisted product $\torus \twisted I$,
which has four boundary faces.  Attached to this boundary are
two new solid tori.  This procedure is illustrated in
Figure~\ref{fig-pluggedthin}.

The fibration results as follows.
Let $\mb$ represent the {\mobius} band, and let the orbifold
$\annref$ be the annulus with one reflector boundary and one regular
boundary component.
The twisted product $\torus \twisted I$ can be represented as a trivial
{\sfslong} over either $\mb$ or $\annref$ (depending upon the placement
of the fibres).\footnote{More generally, the Seifert fibrations
of every $I$-bundle over the torus or Klein bottle are classified by
Amendola and Martelli in Appendix~A of \cite{italian-nor7}.}
The two new tori then
close off the base orbifold and introduce two exceptional fibres.  The
resulting 3-manifold is a {\sfslong} over either
$\rpp$ or $\discref$ with two exceptional fibres.

\subsubsection{Components}

We now describe details of how the separate components of a plugged thin
$I$-bundle are formed.
The original $\torus \twisted I$ is triangulated as a
{\em twisted thin $I$-bundle}, and the two additional tori are
triangulated as {\em layered solid tori}.  We describe each of these
components in turn.

\begin{defn}[Twisted Thin $I$-Bundle] \label{d-twthinibundle}
    A {\em twisted thin $I$-bundle} over the torus
    is a triangulation of the twisted product $\torus \twisted I$
    formed as follows.

    Consider the interval $I = [0,1]$.
    The twisted product $\torus \twisted I$ is naturally foliated by
    surfaces $\torus \times \{x, 1-x\}$ for $0 \leq x \leq \frac12$.
    For all $x \neq \frac12$, this surface is a double cover of the
    torus $\torus \times \{\frac12\}$.

    As in Definition~\ref{d-unthinibundle}, we require that
    this foliation decomposes each individual tetrahedron into
    either triangles or
    quadrilaterals as illustrated in Figure~\ref{fig-ibundletet2}.
    Once more note that every vertex lies on the boundary $F \times \{0,1\}$,
    as does the upper face in the triangular case and the upper and
    lower edges in the quadrilateral case.

    Again we observe that the surface $F \times \{\frac12\}$ meets every
    tetrahedron in precisely one triangle or quadrilateral.  This
    surface is referred to as the {\em central torus} of the $I$-bundle.

    \begin{figure}[htb]
    \centerline{\includegraphics[scale=0.6]{ibundletet.eps}}
    \fcaption{Decomposing a tetrahedron into triangles or quadrilaterals}
    \label{fig-ibundletet2}
    \end{figure}
\end{defn}

\begin{figure}[htb]
\psfrag{Prod}{{\small \twolines{c}{The twisted}{product $\torus \twisted I$}}}
\psfrag{Bdry}{{\small \twolines{c}{The torus}{boundary}}}
\psfrag{Central}{{\small The central torus}}
\psfrag{A}{{\tiny $A$}} \psfrag{B}{{\tiny $B$}} \psfrag{C}{{\tiny $C$}}
\psfrag{D}{{\tiny $D$}} \psfrag{E}{{\tiny $E$}} \psfrag{F}{{\tiny $F$}}
\psfrag{G}{{\tiny $G$}} \psfrag{H}{{\tiny $H$}} \psfrag{J}{{\tiny $J$}}
\centerline{\includegraphics[scale=0.6]{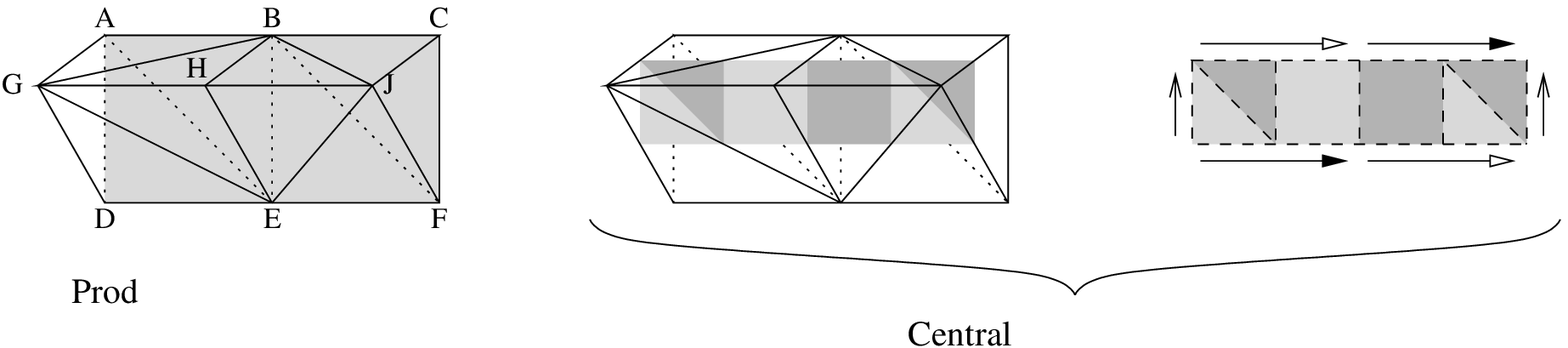}}
\fcaption{An example of a twisted thin $I$-bundle}
\label{fig-ibundle-twt61}
\end{figure}

An example of a twisted thin $I$-bundle over the torus is shown in
Figure~\ref{fig-ibundle-twt61}.  Here we have six tetrahedra arranged
into a long triangular prism, whose four back faces form the boundary
torus (as shaded in the first diagram).
The left and right triangles are identified directly (so that
\polytri{ADG} is identified with \polytri{CFJ}).
The upper and lower
rectangles are identified with a twist and a translation, so that
\polysq{ABHG} and \polysq{HJFE} are identified and
\polysq{GHED} and \polysq{BCJH} are identified.

The central torus $\torus \times \{\frac12\}$ is shaded in the
second diagram of the sequence, and in the third diagram it is made
clear how the six tetrahedra divide this torus into four triangles
and two quadrilaterals.  The arrows on this final diagram indicate which
edges of the central torus are identified with which others.

As with the untwisted thin $I$-bundles of the previous section,
it should be noted that the decomposition of the
central torus into triangles and quadrilaterals
provides enough information to completely reconstruct the thin $I$-bundle.

The second component that appears in a plugged thin $I$-bundle is the
{\em layered solid torus}.  Layered solid tori are well understood, and have
been discussed by Jaco and Rubinstein \cite{0-efficiency,layeredlensspaces}
as well as by Matveev, Martelli and Petronio in the context of special spines
\cite{italianfamilies,matveev6}.
In the context of a census of triangulations they are
discussed and parameterised thoroughly in \cite{burton-nor7}.
We omit the details here.

\begin{figure}[htb]
\psfrag{Torus}{{\small \twolines{c}{Solid}{torus}}}
\centerline{\includegraphics[scale=0.7]{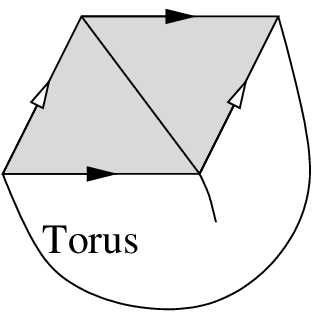}}
\fcaption{The boundary of a layered solid torus}
\label{fig-lst}
\end{figure}

For this overview it suffices to know the following.
A layered solid torus is a triangulation of a solid torus
containing one vertex, two boundary faces and three boundary edges,
as illustrated in Figure~\ref{fig-lst}.
Moreover, it is constructed with the explicit aim of making its
three boundary edges follow some particular curves along the boundary
torus.  There are infinitely many different layered tori,
corresponding to infinitely many different choices of boundary curves.

\begin{figure}[htb]
\psfrag{F}{{\small $F$}}
\psfrag{F'}{{\small $F'$}}
\centerline{\includegraphics[scale=0.7]{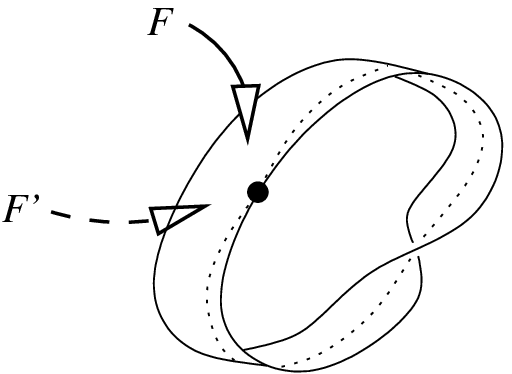}}
\fcaption{A degenerate layered solid torus}
\label{fig-mobius}
\end{figure}

It is useful to consider the {\mobius} band as a
{\em degenerate} layered solid torus with zero tetrahedra.  That is, the
{\mobius} band formed from a single triangle
can be thickened slightly to create a solid torus with two
boundary faces $F$ and $F'$, as illustrated in Figure~\ref{fig-mobius}.

We are now ready to define a plugged thin $I$-bundle.

\begin{defn}[Plugged Thin $I$-Bundle] \label{d-pluggedthin}
    A {\em plugged thin $I$-bundle} is a triangulation constructed as
    follows.  Begin with a twisted thin $I$-bundle over the torus.
    This twisted thin $I$-bundle must have precisely six tetrahedra and
    four boundary faces.  Furthermore, these boundary faces must form
    one of the two configurations shown in Figure~\ref{fig-torus4}.
    We refer to these configurations as the {\em allowable torus
    boundaries}.

    \begin{figure}[htb]
    \psfrag{A}{{\small $A$}} \psfrag{B}{{\small $B$}} \psfrag{C}{{\small $C$}}
    \psfrag{D}{{\small $D$}} \psfrag{E}{{\small $E$}} \psfrag{F}{{\small $F$}}
    \centerline{\includegraphics[scale=0.7]{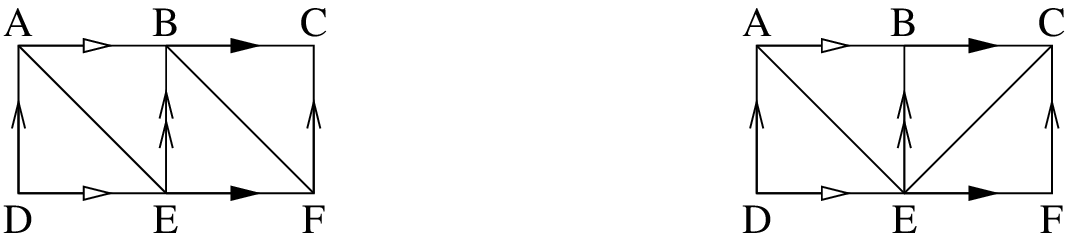}}
    \fcaption{The two allowable torus boundaries}
    \label{fig-torus4}
    \end{figure}

    Observe that these boundary faces can be split into two annuli
    (the left annulus \poly{ABED} and the right annulus \poly{BCFE},
    with edge \polyedge{BE} distinct from edges \polyedge{AD} and
    \polyedge{CF}).  To
    each of these annuli, attach a layered solid torus.  These tori must
    be attached so that edges \polyedge{AD}, \polyedge{BE} and \polyedge{CF}
    are identified, and each annulus \poly{ABED} and \poly{BCFE} becomes a
    torus instead.  This is illustrated for the first boundary
    configuration in Figure~\ref{fig-attachtori}.

    \begin{figure}[htb]
    \centerline{\includegraphics[scale=0.7]{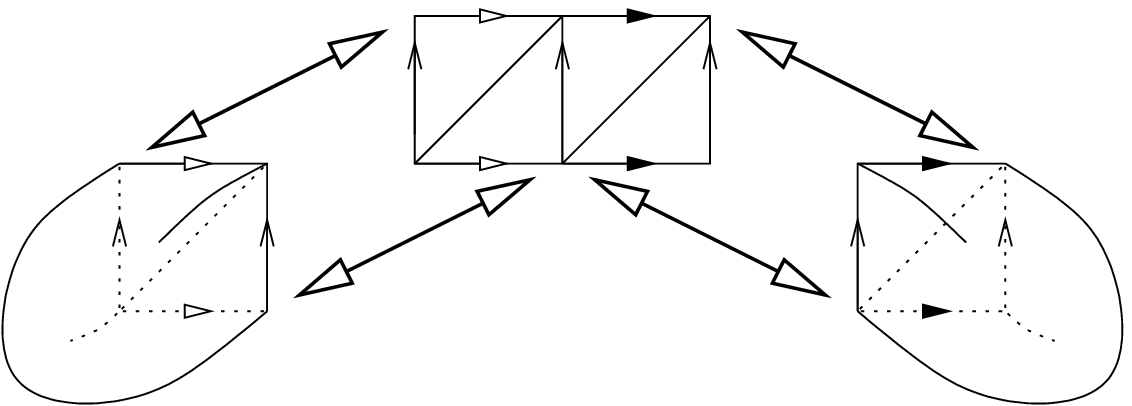}}
    \fcaption{Attaching two layered solid tori to the boundary}
    \label{fig-attachtori}
    \end{figure}

    Note that either layered solid torus may be degenerate.  In this case
    a one-face {\mobius} band is inserted, and the two faces of the
    corresponding boundary annulus are joined to each side of this
    {\mobius} band.  Since the {\mobius} band has no thickness, the
    result is that the two faces of the boundary annulus become joined
    to each other.  An example of this is illustrated in
    Figure~\ref{fig-insertmobius}, where faces \polytri{CDA} and
    \polytri{DAB} become identified.

    \begin{figure}[htb]
    \psfrag{A}{{\small $A$}} \psfrag{B}{{\small $B$}}
    \psfrag{C}{{\small $C$}} \psfrag{D}{{\small $D$}}
    \centerline{\includegraphics[scale=0.6]{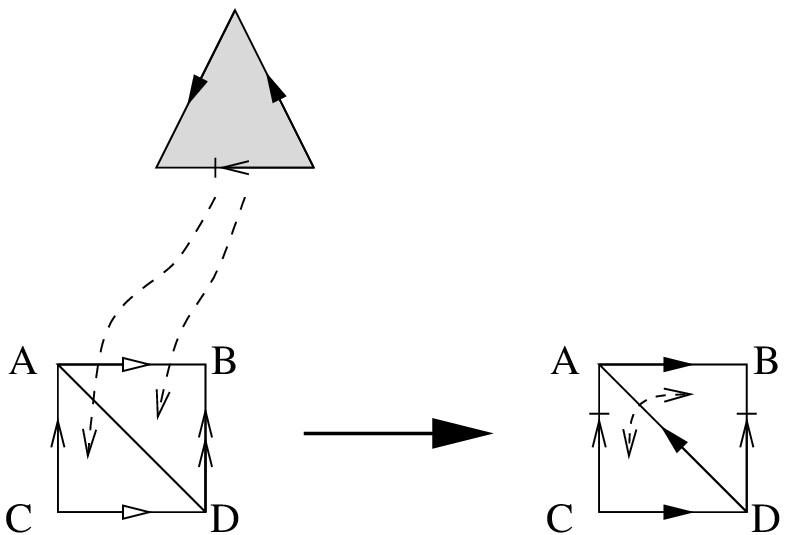}}
    \fcaption{Attaching a degenerate layered solid torus to an annulus}
    \label{fig-insertmobius}
    \end{figure}
\end{defn}

Again it may help to refer to Figure~\ref{fig-pluggedthin} for an
overview of this construction.  The underlying 3-manifold can then be
determined as follows.

\begin{figure}[htb]
\centerline{\includegraphics[scale=0.7]{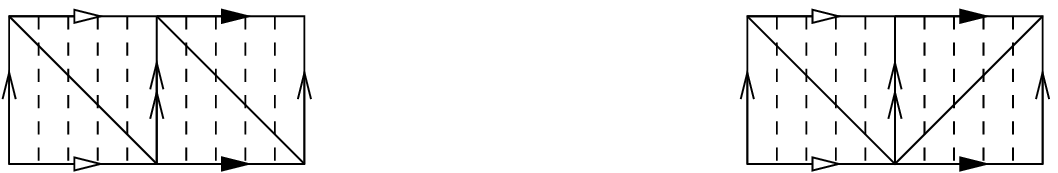}}
\fcaption{Fibres in the two allowable torus boundaries}
\label{fig-sfsvert4}
\end{figure}

The four-face
boundary of the twisted $\torus \twisted I$ can be filled with vertical
fibres, as illustrated in Figure~\ref{fig-sfsvert4}.  As shown by
Amendola and Martelli \cite{italian-nor7}, this extends to a Seifert
fibration
of $\torus \twisted I$ as a trivial {\sfslong} over either $\mb$ or
$\annref$.  In the other direction, this extends to a
Seifert fibration of each layered solid torus with an
exceptional fibre at its centre (unless the boundary curves for the
layered solid torus are chosen so that the meridinal disc of the torus is
bounded by a fibre or meets each fibre just once).

The result is a {\sfslong} over either $\rpp$ or $\discref$ with two
exceptional fibres.  See \cite{burton-nor7} for a formula that gives
the precise Seifert invariants in terms of the individual parameters of
the thin $I$-bundle and the two layered solid tori.

\subsubsection{Census Triangulations}

The plugged thin $I$-bundles in the census give rise to the
{\sfslong}s
\[ \begin{array}{lll}
\sfs{\rpp}{(2,1)\ (2,1)}, \quad &
\sfs{\rpp}{(2,1)\ (3,1)}, \quad &
\sfs{\rpp}{(2,1)\ (4,1)}, \\
\sfs{\rpp}{(2,1)\ (5,2)}, \quad &
\sfs{\rpp}{(3,1)\ (3,1)}, \quad &
\sfs{\rpp}{(3,1)\ (3,2)}, \\
\sfs{\discref}{(2,1)\ (2,1)}, \quad &
\sfs{\discref}{(2,1)\ (3,1)}, \quad &
\sfs{\discref}{(2,1)\ (4,1)}, \\
\sfs{\discref}{(2,1)\ (5,2)}, \quad &
\sfs{\discref}{(3,1)\ (3,1)}, \quad &
\sfs{\discref}{(3,1)\ (3,2)}.
\end{array} \]

Table~\ref{tab-thin} lists the frequencies of plugged thin $I$-bundles
within the overall census.  Here the large number of triangulations
results from the fact that there are several possible choices for the
initial twisted product $\torus \twisted I$,
as well as the equivalence between some
spaces such as $\sfs{\rpp}{(3,1)\ (3,1)}$ and $\sfs{\rpp}{(3,2)\ (3,2)}$.

\begin{table}[htb]
\tcaption{Frequencies of plugged thin $I$-bundles within the census}
\begin{center} \begin{tabular}{|r|rl|rl|}
    \hline
    \bf Tetrahedra & \multicolumn{2}{c|}{\bf 3-Manifolds} &
    \multicolumn{2}{c|}{\bf Triangulations} \\
    \hline
    6 & 2 & (out of 5) & 4 & (out of 24) \\
    7 & 2 & (out of 3) & 6 & (out of 17) \\
    8 & 8 & (out of 10) & 22 & (out of 59) \\
    \hline
\end{tabular} \end{center}
\label{tab-thin}
\end{table}

%
%

\subsection{Plugged Thick $I$-Bundles} \label{s-tri-thick}

\begin{figure}[htb]
\psfrag{TP}{{\small Twisted product $\torus \twisted I$ (3 or 5 tets)}}
\psfrag{Bdry2}{{\small Bad boundary torus (2 or 4 faces)}}
\psfrag{Plug}{{\small Thickening plug}}
\psfrag{Bdry4}{{\small Good boundary torus (4 faces)}}
\psfrag{LST}{{\small \twolines{ll}{Two layered}{solid tori}}}
\centerline{\includegraphics[scale=0.6]{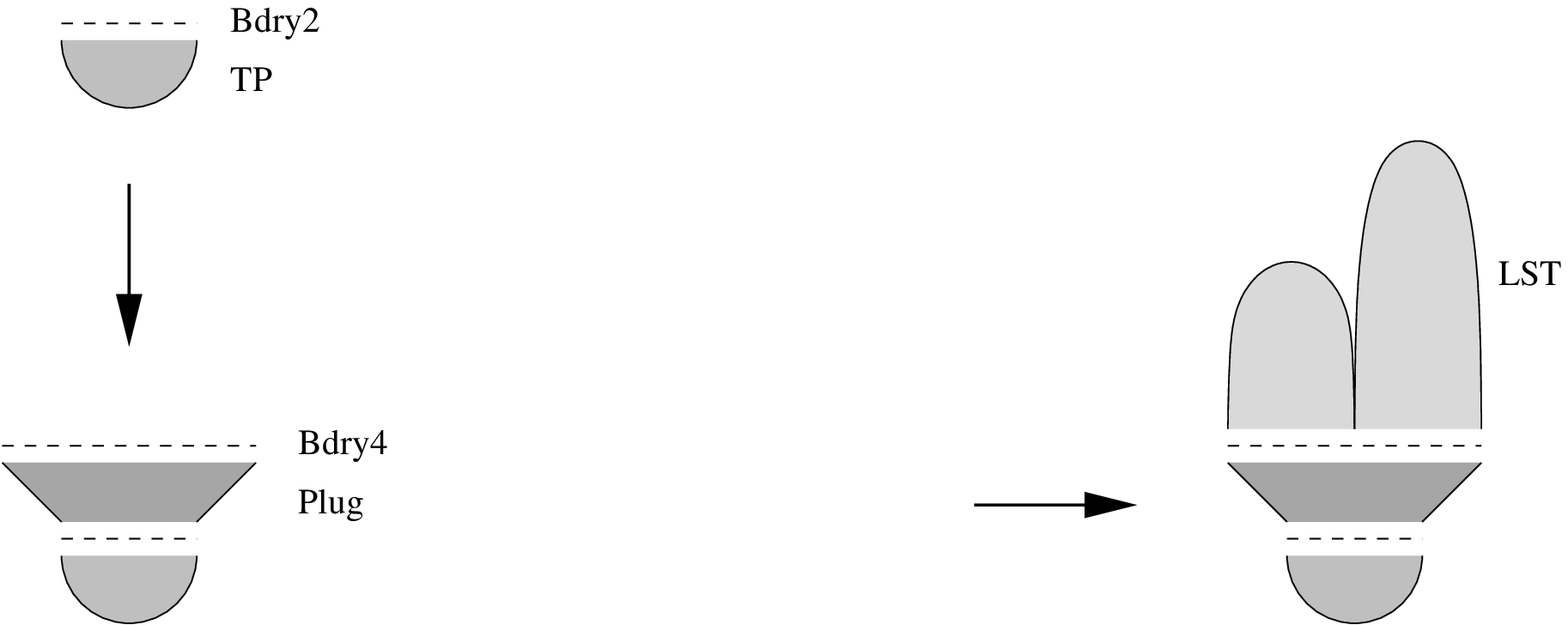}}
\fcaption{Constructing a plugged thick $I$-bundle}
\label{fig-pluggedthick}
\end{figure}

A {\em plugged thick $I$-bundle} is very similar in construction to a plugged
thin $I$-bundle.  The difference is that a smaller twisted thin
$I$-bundle is used, but the resulting torus boundary is not one of the
allowable torus boundaries of Figure~\ref{fig-torus4}.  As a result, some
additional tetrahedra must be added to reconfigure the torus boundary
(thus ``thickening'' the $I$-bundle).  Once this is done,
the two new layered solid tori are attached as before.
This procedure is illustrated in Figure~\ref{fig-pluggedthick}.

There are two different ways in which this construction can be carried
out.
\begin{enumerate}[(i)]
    \item
    {\em We begin with a three-tetrahedron twisted thin $I$-bundle over the
    torus.}  An example is shown in Figure~\ref{fig-ibundle-t31}.  The
    three tetrahedra are arranged into a triangular prism, and the
    boundary torus is formed from the two back faces (as shaded in the
    second diagram).  The left and right faces are identified directly (with
    \polytri{ACE} identified with \polytri{BDF}), and the upper and
    lower squares are identified with a twist and a translation
    (with \polytri{EFA} and \polytri{CDF} identified and with
    \polytri{ABF} and \polytri{EFC} identified).
    As usual, the third and fourth diagrams illustrate the
    the central torus $\torus \times \{\frac12\}$.

    \begin{figure}[htb]
    \psfrag{Prod}{{\small \twolines{c}{The twisted}{product $\torus \twisted I$}}}
    \psfrag{Bdry}{{\small \twolines{c}{The torus}{boundary}}}
    \psfrag{Central}{{\small The central torus}}
    \psfrag{A}{{\tiny $A$}} \psfrag{B}{{\tiny $B$}} \psfrag{C}{{\tiny $C$}}
    \psfrag{D}{{\tiny $D$}} \psfrag{E}{{\tiny $E$}} \psfrag{F}{{\tiny $F$}}
    \centerline{\includegraphics[scale=0.6]{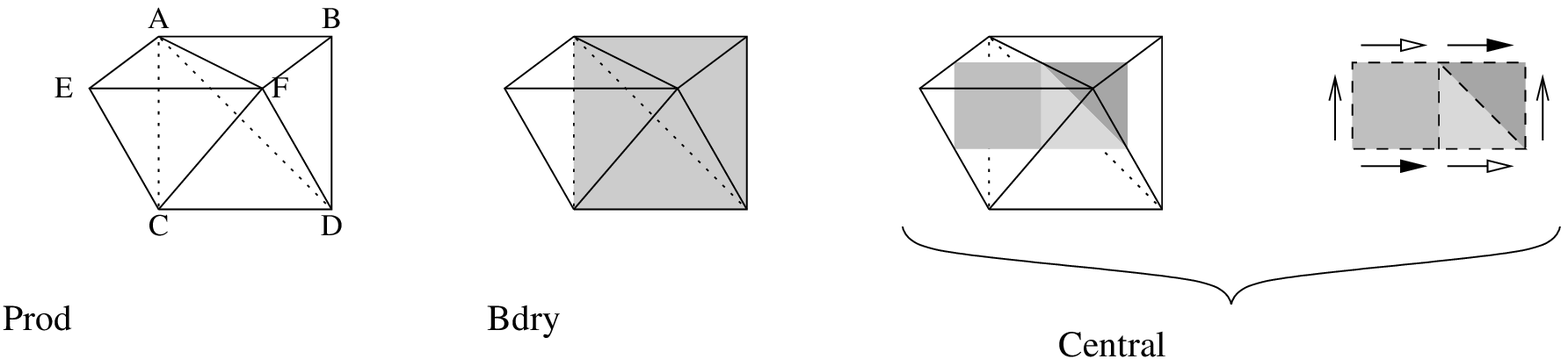}}
    \fcaption{A three-tetrahedron twisted thin $I$-bundle}
    \label{fig-ibundle-t31}
    \end{figure}

    The resulting boundary torus has only two faces, which is clearly
    not an allowable torus boundary (see Definition~\ref{d-pluggedthin}
    and Figure~\ref{fig-torus4}).
    To compensate, we attach a three-tetrahedron thickening plug as
    illustrated in Figure~\ref{fig-plug3}.  The two back faces of this
    plug (shaded in the diagram) are attached to the old two-face boundary
    torus.  The four front faces become a new allowable boundary torus,
    and the upper and lower faces \polytri{ABC} and \polytri{DEF} are
    identified with each other.

    \begin{figure}
    \psfrag{A}{{\tiny $A$}} \psfrag{B}{{\tiny $B$}} \psfrag{C}{{\tiny $C$}}
    \psfrag{D}{{\tiny $D$}} \psfrag{E}{{\tiny $E$}} \psfrag{F}{{\tiny $F$}}
    \centerline{\includegraphics[scale=0.7]{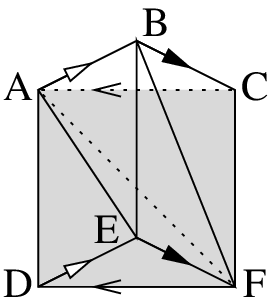}}
    \fcaption{A three-tetrahedron thickening plug}
    \label{fig-plug3}
    \end{figure}

    \item
    {\em We begin with a five-tetrahedron twisted thin $I$-bundle over
    the torus.}  This is illustrated in Figure~\ref{fig-ibundle-t51},
    with the faces of the triangular prism identified as in
    Figure~\ref{fig-ibundle-t31} before.
    Here the torus boundary has four faces, but they are not arranged
    into an allowable torus boundary.  We must therefore reconfigure
    the boundary edges by performing a layering, resulting in a new
    four-face boundary that satisfies our requirements.

    \begin{figure}[htb]
    \psfrag{Prod}{{\small \twolines{c}{The twisted}{product $\torus \twisted I$}}}
    \psfrag{Bdry}{{\small \twolines{c}{The torus}{boundary}}}
    \psfrag{Central}{{\small The central torus}}
    \centerline{\includegraphics[scale=0.6]{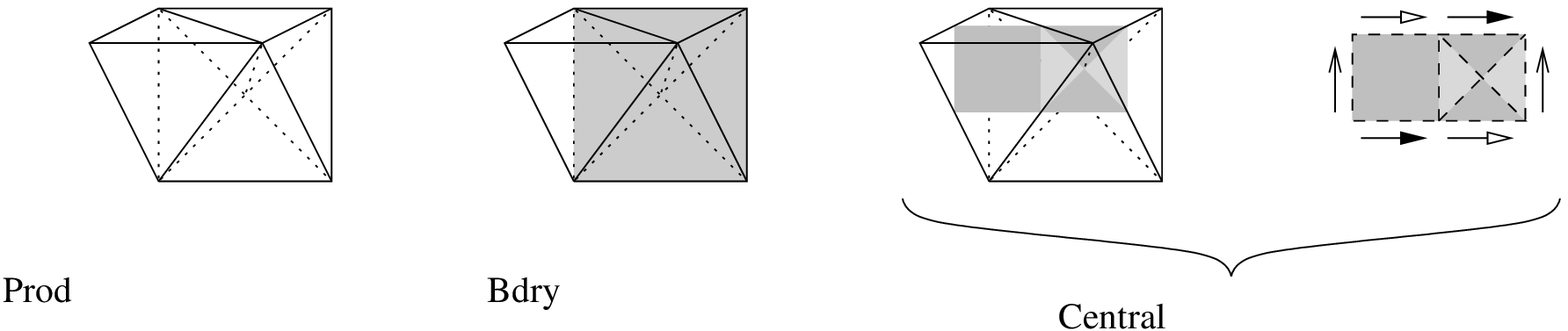}}
    \fcaption{A five-tetrahedron twisted thin $I$-bundle}
    \label{fig-ibundle-t51}
    \end{figure}
\end{enumerate}

We can formalise this into the following definition.  For a complete
enumeration of the different $I$-bundles and thickening plugs that can be
used, the reader is referred to \cite{burton-nor7}.

\begin{defn}[Plugged Thick $I$-Bundle] \label{d-pluggedthick}
    A {\em plugged thick $I$-bundle} is a triangulation formed as
    follows.  We begin with a twisted thin $I$-bundle over the torus,
    which either has (i) three tetrahedra and two boundary faces, or
    (ii) five tetrahedra and four boundary faces.  We then convert the
    torus boundary into one of the allowable torus boundaries described
    in Definition~\ref{d-pluggedthin} (see Figure~\ref{fig-torus4}),
    either by (i) inserting a three-tetrahedron
    thickening plug as described above, or (ii) layering a new
    tetrahedron onto the torus boundary.

    We require that
    the resulting structure is a six-tetrahedron triangulation
    of the twisted product $\torus \twisted I$ with an allowable torus
    boundary.  We finish the construction by attaching two layered solid
    tori exactly as described in Definition~\ref{d-pluggedthin}.
\end{defn}

Since we are producing
triangulations of the twisted product $\torus \twisted I$ with the same
allowable boundary tori used for plugged thin $I$-bundles,
it follows that we should obtain the same underlying
3-manifolds.  Specifically, we obtain
{\sfslong}s over either $\rpp$ or $\discref$ with two
exceptional fibres.  Again a precise formula appears in \cite{burton-nor7}
to calculate the exact Seifert invariants in terms of the individual
parameters of the triangulation.

\subsubsection{Census Triangulations}

The plugged thick $I$-bundles in the census give rise to the spaces
\[ \begin{array}{lll}
\sfs{\rpp}{(2,1)\ (2,1)}, \quad &
\sfs{\rpp}{(2,1)\ (3,1)}, \quad &
\sfs{\rpp}{(2,1)\ (4,1)}, \\
\sfs{\rpp}{(2,1)\ (5,2)}, \quad &
\sfs{\rpp}{(3,1)\ (3,1)}, \quad &
\sfs{\rpp}{(3,1)\ (3,2)}, \\
\sfs{\discref}{(2,1)\ (2,1)}, \quad &
\sfs{\discref}{(2,1)\ (3,1)}, \quad &
\sfs{\discref}{(2,1)\ (4,1)}, \\
\sfs{\discref}{(2,1)\ (5,2)}, \quad &
\sfs{\discref}{(3,1)\ (3,1)}, \quad &
\sfs{\discref}{(3,1)\ (3,2)}.
\end{array} \]
As expected from the similarity in construction, these are
exactly the same 12 spaces that the plugged thin $I$-bundles produce
(though none of the specific triangulations are the
same).  Table~\ref{tab-thick} lists the frequencies of plugged
thick $I$-bundles within the overall census.

\begin{table}[htb]
\tcaption{Frequencies of plugged thick $I$-bundles within the census}
\begin{center} \begin{tabular}{|r|rl|rl|}
    \hline
    \bf Tetrahedra & \multicolumn{2}{c|}{\bf 3-Manifolds} &
    \multicolumn{2}{c|}{\bf Triangulations} \\
    \hline
    6 & 2 & (out of 5) & 4 & (out of 24) \\
    7 & 2 & (out of 3) & 7 & (out of 17) \\
    8 & 8 & (out of 10) & 25 & (out of 59) \\
    \hline
\end{tabular} \end{center}
\label{tab-thick}
\end{table}

\section{Observations and Conjectures} \label{s-obs}

In this final section, we pull together observations from the census
and form conjectures based upon these observations.  Following this we
discuss the future of the non-orientable census, including what we might
expect to see when the census is extended to higher numbers of tetrahedra.

As explained in the introduction, there is an extremely heavy
computational load in creating a census such as this.  Each new level of
the census (measured by number of tetrahedra, or equivalently by the
complexity of Matveev \cite{matveev-complexity}) is an order of magnitude
more difficult to construct than the last.  At the time of writing the
nine-tetrahedron non-orientable census is under construction, with a
healthy body of partial results already obtained.  Each of the
conjectures below is consistent with these partial results.
We return specifically to the nine-tetrahedron census in
Section~\ref{s-obs-future}.

\subsection{Minimal Triangulations} \label{s-obs-minimal}

Our first observation relates to the combinatorial structure of
non-orientable minimal triangulations.  Recall from the introduction
that very few sufficient conditions are known for minimal
triangulations.  Conjectures have been made for various classes of
3-manifolds, but such conjectures are notoriously difficult to prove.

Matveev \cite{matveev6} and Martelli and Petronio \cite{italianfamilies}
have made a variety of well-grounded conjectures about the smallest
number of tetrahedra required for various classes of orientable
3-manifolds.
Here we form conjectures of this type in the
non-orientable case.  Moreover, based upon the growing body of
experimental evidence, we push further and make conjectures regarding the
construction of all minimal triangulations of various classes of
non-orientable 3-manifolds.

Section~\ref{s-tri} introduces
three families of triangulations: (i)~layered surface bundles,
(ii)~plugged thin $I$-bundles and (iii)~plugged thick $I$-bundles.
It is easy enough to see that these families produce
(i)~torus or Klein bottle bundles over the circle and
(ii,iii)~{\sfslong}s over $\rpp$ or $\discref$ with two exceptional
fibres.  What is less obvious is that {\em every} minimal triangulation
of such a 3-manifold should belong to one of the three families listed
above.

In fact the evidence does support this suggestion, with
the exception of the four flat manifolds at the lowest (six-tetrahedron)
level of the census.  Although most 3-manifolds have many different
minimal triangulations (up to 10 in some cases), these
minimal triangulations all belong to
the three families above.  The partial results for the nine-tetrahedron
census also support this hypothesis, even though a much wider variety of
triangulations is found at this level (as discussed below in
Section~\ref{s-obs-future}).
We are therefore led to make the following conjectures.

\begin{conjecture} \label{cj-lsb}
    Let $M$ be a torus bundle over the circle that is not one of the
    flat manifolds $\torus \times I / \homtwo{0}{1}{1}{0}$ or
    $\torus \times I / \homtwo{1}{0}{0}{-1}$.  Then every minimal
    triangulation of $M$ is a layered torus bundle, as described by
    Definition~\ref{d-lsb}.

    Moreover, at least one minimal triangulation of $M$ has at its core
    the six-tetrahedron product $\torus \times I$ illustrated in
    Figure~\ref{fig-lsb-t62}.  In other words, this six-tetrahedron
    $\torus \times I$ may be used as a starting point for constructing a
    minimal triangulation of $M$.
\end{conjecture}

\begin{conjecture} \label{cj-sfs}
    Let $M$ be a {\sfslong} over either $\rpp$ or $\discref$ with
    precisely two exceptional fibres.  Moreover, suppose that $M$ is not
    one of the flat manifolds
    $\sfs{\rpp}{(2,1)\ (2,1)}$ or $\sfs{\discref}{(2,1)\ (2,1)}$.
    Then every minimal triangulation of $M$ is either a plugged thin
    $I$-bundle or a plugged thick $I$-bundle, as described by
    Definitions~\ref{d-pluggedthin} and~\ref{d-pluggedthick}.
\end{conjecture}

Note that if these conjectures are true, the number of tetrahedra in
such a minimal triangulation is straightforward to calculate.  The
number of layerings require to obtain a particular set of boundary
curves is well described by Martelli and Petronio \cite{italianfamilies}
(though in the equivalent language of special spines).  Similar
calculations in the language of triangulations and layered solid tori
have been described by Jaco and Rubinstein in a variety of informal
contexts.

What remains then is to calculate the number of tetrahedra that are not
involved in layerings.  For Conjecture~\ref{cj-lsb} we can assume this
to be the six-tetrahedron $\torus \times I$ of Figure~\ref{fig-lsb-t62},
and for Conjecture~\ref{cj-sfs} we can simply count
the six additional tetrahedra involved in the twisted product
$\torus \twisted I$ to which our layered solid tori are attached.

\subsection{Central Surfaces}

Our next observation is regarding embedded surfaces within
non-orientable 3-manifolds.  Recall that all three families of
triangulations described in Section~\ref{s-tri} begin with a thin
$I$-bundle.  The central surface of this thin $I$-bundle is an embedded
surface meeting each tetrahedron of the thin $I$-bundle in either a
single quadrilateral or a single triangle.

\begin{figure}
\centerline{\includegraphics[scale=0.6]{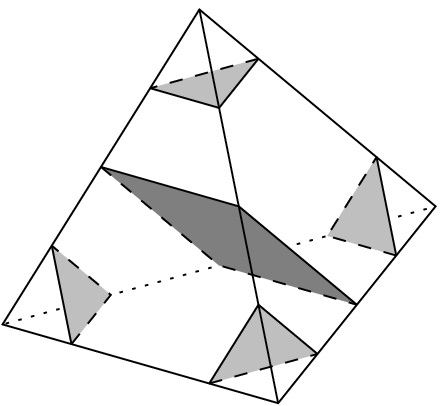}}
\fcaption{Normal discs within a tetrahedron}
\label{fig-normaldiscs}
\end{figure}

This is reminiscent of the theory of normal surfaces.  Normal surfaces,
first introduced by Kneser \cite{kneser-normal} and subsequently developed
by Haken \cite{haken-knot,haken-homeomorphism}, play a powerful role in
algorithms in 3-manifold topology.  A {\em normal surface} within a
triangulation meets each tetrahedron in one or more {\em normal
discs}, which are either triangles separating one vertex from the other
three or quadrilaterals separating two vertices from the other two.  A
variety of normal discs can be seen in Figure~\ref{fig-normaldiscs}.

It follows then that the central surface of a thin $I$-bundle is a
special type of normal surface, namely one that meets each tetrahedron
of the thin $I$-bundle in one and only one normal disc.  This leads us
to make the following more general definition.

\begin{defn}[Central Normal Surface]
    Let $T$ be a 3-manifold triangulation, and let $N$ be an embedded
    normal surface in $T$.  We refer to $N$ as a {\em central normal
    surface} if and only if $N$ meets each tetrahedron in at most one
    normal disc (i.e., one triangle, one quadrilateral or nothing).
\end{defn}

Note that we have replaced ``one and only one'' with ``at most one'',
since the central surface of a thin $I$-bundle does not meet the
tetrahedra involved in the other parts of the triangulation (such as the
layered solid tori in a plugged thin $I$-bundle).

It can be seen that every triangulation in this census contains a
central normal surface.  Indeed this is to be expected, given that a
thin $I$-bundle appears at the core of every family described in
Section~\ref{s-tri}.  However, every triangulation obtained thus far
in the nine-tetrahedron census also contains a central normal surface,
even though a much greater variety of 3-manifolds and triangulations are
obtained.

To account for the possibility that the central normal surfaces are
specific features of {\sfslong}s or more generally of graph manifolds,
a number of hyperbolic manifolds have also been examined.  As part of a
census of closed hyperbolic 3-manifolds \cite{closedhypcensus},
Hodgson and Weeks list ten candidates for the smallest volume closed
non-orientable hyperbolic 3-manifolds.  An extended version of this
table is shipped with {\snappea} \cite{snappea}, giving a total of 18 closed
non-orientable hyperbolic 3-manifolds with volumes ranging from
$2.02988321$ to $3.85787307$.

Additional experimentation on triangulations of these hyperbolic
3-manifolds yields a similar result, namely that each contains a
central normal surface.  We therefore make the following conjecture.

\begin{conjecture} \label{cj-central}
    Every minimal triangulation of a closed non-orientable
    {\ppirr} 3-manifold contains a central normal surface.
\end{conjecture}

It is worth noting that this conjecture does not hold in the orientable
case.  Of the 191 closed orientable minimal irreducible triangulations
with $\leq 6$ tetrahedra\footnote{These triangulations are enumerated by
Matveev in \cite{matveev6} in the equivalent language of special spines.},
only 118 have a central normal surface.
It is known that every
closed non-orientable 3-manifold contains an embedded incompressible
surface \cite{hempel}, which may help explain the presence of central
normal surfaces in the non-orientable case.

It may be that Conjecture~\ref{cj-central} is easier to prove than
Conjectures~\ref{cj-lsb} and~\ref{cj-sfs}, since existence theorems for
normal surfaces are generally easier to come by.  Moreover, once proven,
Conjecture~\ref{cj-central} may well be a starting point for the proofs
of the other conjectures --- as with the central surface of a thin
$I$-bundle, a central normal surface can be used to reconstruct the
portion of the triangulation that surrounds it, which may then lead
to new structural results.

\subsection{Future Directions} \label{s-obs-future}

As mentioned in the introduction, the triangulations and 3-manifolds
seen in the eight-tetrahedron census offer little variety
beyond what has already been seen in the $\leq 7$-tetra\-hedron
census.  The primary advantage of the eight-tetrahedron census has been
the larger body of data (an additional 59 triangulations and 10 distinct
3-manifolds) that has supported the formulation of the
conjectures above.

Clearly a greater variety must appear in the census at some point, since
there are far more non-orientable 3-manifolds than those described by
the families of Section~\ref{s-tri}.  The question is how much
larger the census must become before we begin to see them.

Fortunately the answer is not much larger at all.  As discussed at the
beginning of Section~\ref{s-obs}, the nine-tetrahedron census is
currently under construction and a significant body of partial results
are already available.  In addition to the 3-manifolds already described
(torus bundles over the circle and {\sfslong}s over $\rpp$ or $\discref$
with two exceptional fibres), the nine-tetrahedron results include the
following:
\begin{itemize}
    \item {\em {\sfslong}s over $\rpp$ and $\discref$ with three exceptional
    fibres.}  In particular, the spaces $\sfs{\rpp}{(2,1)\ (2,1)\ (2,1)}$ and
    $\sfs{\discref}{(2,1)\ (2,1)\ (2,1)}$ are found.

    \item {\em {\sfslong}s over several other base orbifolds with one
    exceptional fibre.}  The base orbifolds include the torus, the Klein
    bottle, the annulus with two reflector boundaries, and the {\mobius}
    strip with one reflector boundary.  In each case a single $(2,1)$
    exceptional fibre is found.

    \item {\em Manifolds with non-trivial JSJ composition.}  In
    particular, a number of spaces are found that begin with a
    {\sfslong} over the annulus with a single $(2,1)$ fibre, followed
    by a non-trivial identification of the two torus boundaries.
\end{itemize}

It is therefore hoped that, once completed, the nine-tetrahedron
census can offer richer insights into the structures of non-orientable
minimal triangulations than what we have seen to date.

Moving beyond the nine-tetrahedron census, one might ask how much
further we must go before we move away from graph manifolds.  It is
noted by Martelli and Petronio \cite{italianfamilies} that the first
hyperbolic manifolds to appear in the orientable census are those of
smallest known volume (first seen at nine tetrahedra).  It is reasonable
to expect the same of the non-orientable census; the smallest volume
non-orientable hyperbolic manifold described by Hodgson and Weeks
\cite{closedhypcensus} can be triangulated with 11 tetrahedra, though
neither the minimality of the volume nor the minimality of
the triangulation have been proven.

Finally it must be noted that any extension of the census will require
new improvements in the algorithm --- the worse-than-exponential growth
of the search space means that increased computing power is not enough.
Ideally such improvements would involve a blend of topological results
(such as those seen in \cite{burton-facegraphs}) and pure algorithmic
optimisations.  The expected yield from higher levels of the census is a
great incentive, and so work on the enumeration algorithm is continuing.

\bibliographystyle{amsplain}
\bibliography{pure-20050910}

\newcommand{\noopsort}[1]{}
\providecommand{\bysame}{\leavevmode\hbox to3em{\hrulefill}\thinspace}
\providecommand{\MR}{\relax\ifhmode\unskip\space\fi MR }
\providecommand{\MRhref}[2]{%
  \href{http://www.ams.org/mathscinet-getitem?mr=#1}{#2}
}
\providecommand{\href}[2]{#2}
\begin{thebibliography}{10}

\bibitem{italian-nor6}
Gennaro Amendola and Bruno Martelli, \emph{Non-orientable 3-manifolds of small
  complexity}, Topology Appl. \textbf{133} (2003), no.~2, 157--178.

\bibitem{italian-nor7}
\bysame, \emph{Non-orientable 3-manifolds of complexity up to 7}, Topology
  Appl. \textbf{150} (2005), no.~1-3, 179--195.

\bibitem{regina}
Benjamin~A. Burton, \emph{Regina: Normal surface and 3-manifold topology
  software}, {\tt http://\allowbreak regina.\allowbreak sourceforge.\allowbreak
  net/}, 1999--2005.

\bibitem{burton-nor7}
\bysame, \emph{Structures of small closed non-orientable 3-manifold
  triangulations}, Preprint, {\tt math.\allowbreak GT/\allowbreak 0311113},
  November 2003.

\bibitem{burton-facegraphs}
\bysame, \emph{Face pairing graphs and 3-manifold enumeration}, J. Knot Theory
  Ramifications \textbf{13} (2004), no.~8, 1057--1101.

\bibitem{burton-regina}
\bysame, \emph{Introducing {R}egina, the 3-manifold topology software},
  Experiment. Math. \textbf{13} (2004), no.~3, 267--272.

\bibitem{cuspedcensus}
Patrick~J. Callahan, Martin~V. Hildebrand, and Jeffrey~R. Weeks, \emph{A census
  of cusped hyperbolic 3-manifolds}, Math. Comp. \textbf{68} (1999), no.~225,
  321--332.

\bibitem{italian-heegaard}
Roberto Frigerio, Bruno Martelli, and Carlo Petronio, \emph{Complexity and
  {H}eegaard genus of an infinite class of compact 3-manifolds}, Pacific J.
  Math. \textbf{210} (2003), no.~2, 283--297.

\bibitem{haken-knot}
Wolfgang Haken, \emph{Theorie der {N}ormalfl\"achen}, Acta Math. \textbf{105}
  (1961), 245--375.

\bibitem{haken-homeomorphism}
\bysame, \emph{\"{U}ber das {H}om\"oomorphieproblem der 3-{M}annigfaltigkeiten.
  {I}}, Math. Z. \textbf{80} (1962), 89--120.

\bibitem{hempel}
John Hempel, \emph{3-manifolds}, Annals of Mathematics Studies, no.~86,
  Princeton University Press, Princeton, NJ, 1976.

\bibitem{closedhypcensus}
Craig~D. Hodgson and Jeffrey~R. Weeks, \emph{Symmetries, isometries and length
  spectra of closed hyperbolic three-manifolds}, Experiment. Math. \textbf{3}
  (1994), no.~4, 261--274.

\bibitem{0-efficiency}
William Jaco and J.~Hyam Rubinstein, \emph{0-efficient triangulations of
  3-manifolds}, J. Differential Geom. \textbf{65} (2003), no.~1, 61--168.

\bibitem{layeredlensspaces}
\bysame, \emph{Layered triangulations of lens spaces}, In preparation, 2003.

\bibitem{kneser-normal}
Hellmuth Kneser, \emph{Geschlossene {F}l\"achen in dreidimensionalen
  {M}annigfaltigkeiten}, Jahresbericht der Deut. Math. Verein. \textbf{38}
  (1929), 248--260.

\bibitem{italian10}
Bruno Martelli, \emph{Complexity of 3-manifolds}, Preprint, {\tt
  math.\allowbreak GT/\allowbreak 0405250}, January 2005.

\bibitem{italian9}
Bruno Martelli and Carlo Petronio, \emph{Three-manifolds having complexity at
  most 9}, Experiment. Math. \textbf{10} (2001), no.~2, 207--236.

\bibitem{italian-decomp}
\bysame, \emph{A new decomposition theorem for 3-manifolds}, Illinois J. Math.
  \textbf{46} (2002), 755--780.

\bibitem{italianfamilies}
\bysame, \emph{Complexity of geometric three-manifolds}, Geom. Dedicata
  \textbf{108} (2004), no.~1, 15--69.

\bibitem{matveev-complexity}
Sergei~V. Matveev, \emph{Complexity theory of three-dimensional manifolds},
  Acta Appl. Math. \textbf{19} (1990), no.~2, 101--130.

\bibitem{matveev6}
\bysame, \emph{Tables of 3-manifolds up to complexity 6}, Max-Planck-Institut
  f\"{u}r Mathematik Preprint Series (1998), no.~67, available from {\tt
  http://www.\allowbreak mpim-bonn.\allowbreak mpg.\allowbreak de/\allowbreak
  html/\allowbreak pre\-prints/\allowbreak preprints.html}.

\bibitem{snappea}
Jeffrey~R. Weeks, \emph{Snap{P}ea ({H}yperbolic 3-manifold software)}, {\tt
  http://www.\allowbreak northnet.\allowbreak org/\allowbreak weeks/\allowbreak
  index/\allowbreak SnapPea.html}, 1991--2000.

\end{thebibliography}

\section*{Appendix}

For convenience we include a list of all 100 triangulations from
the census, named according to the precise parameterisations described
in \cite{burton-nor7}.  This allows the reader to cross-reference
triangulations and constructions between these two papers.  In summary
we have the following.
\begin{itemize}
    \item Triangulations $H_{\twt\ldots}$ are plugged thin $I$-bundles, and
    triangulations $K_{\twt\ldots}$ are plugged thick $I$-bundles.

    \item Triangulations $B_{T\ldots}$ are layered torus bundles, and
    triangulations $B_{K\ldots}$ are layered Klein bottle bundles.
    Note that there is one triangulation of
    $\torus \times I / \homtwo{0}{1}{1}{0}$ that can be expressed in
    both forms.

    \item Triangulations $E_{6,1}$, $E_{6,2}$ and $E_{6,3}$ are
    described in \cite{burton-nor7} as exceptional triangulations.
    However, both $E_{6,1}$ and $E_{6,2}$ are also layered Klein bottle
    bundles whose central Klein bottles were not originally included in
    the parameterisation of \cite{burton-nor7}.

    \item More specifically,
    triangulations $B_{T_n\ldots}$ and $B_{K_n\ldots}$ are
    constructed from thin $I$-bundles containing precisely $n$ tetrahedra.
    The four triangulations named $B_{T_8\ldots}$ are not
    included in the parameterisation of \cite{burton-nor7}, since
    eight-tetrahedron thin $I$-bundles were not covered.
\end{itemize}
For full details, including a precise explanation of the parameterisation
system, the reader is referred to \cite{burton-nor7}.

\renewcommand{\arraystretch}{1.5}
\begin{table}[htb]
\tcaption{All 18 distinct 3-manifolds and their 100 minimal triangulations}
\small
\begin{center} \begin{tabular}{|c|l|l|}
    \hline
    $\Delta$ & \bf 3-Manifold & \bf Triangulations \\
    \hline \hline
6 & $\torus \times I / \homtwo{1}{1}{1}{0}$ &
    $B_{T_6^2 | -1,1|1,0}$ \\

  \cline{2-3}
  & $\torus \times I / \homtwo{0}{1}{1}{0}$ &
    $B_{T_6^1 | -1,0|-1,1}$,
    $B_{T_6^1 | 0,-1|-1,0}$,
    $B_{T_6^1 | 0,1|1,0} = B_{K_6^2 | 0,-1|-1,0}$, \\
  & &
    $B_{T_6^1 | 1,0|1,-1}$,
    $B_{K_6^1 | 0,-1|-1,0}$,
    $E_{6,3}$ \\

  \cline{2-3}
  & $\torus \times I / \homtwo{1}{0}{0}{-1}$ &
    $B_{T_6^2 | 1,0|0,-1}$,
    $B_{K_6^1 | 1,0|0,1}$,
    $B_{K_6^2 | 1,0|0,1}$ \\

  \cline{2-3}
  & $\sfs{\rpp}{(2,1)\ (2,1)}$ &
    $B_{K_6^1 | 0,1|1,0}$,
    $B_{K_6^2 | 0,1|1,0}$,
    $H_{\twt_6^1}$,
    $H_{\twt_6^2}$,
    $H_{\twt_6^3}$, \\
  & &
    $K_{\twt_5^1}$,
    $K_{\twt_5^2}$,
    $K_{\twt_5^3}$,
    $E_{6,2}$ \\

  \cline{2-3}
  & $\sfs{\discref}{(2,1)\ (2,1)}$ &
    $B_{K_6^1 | -1,0|0,-1}$,
    $B_{K_6^2 | -1,0|0,-1}$,
    $H_{\twt_6^4}$,
    $K_{\twt_5^4}$,
    $E_{6,1}$ \\
  \hline \hline

7 & $\torus \times I / \homtwo{2}{1}{1}{0}$ &
    $B_{T_6^2 | -1,1 | 2,-1}$,
    $B_{T_6^2 | 0,-1 | -1,2}$,
    $B_{T_7 | -1,-1 | -1,0}$,
    $B_{T_7 | 1,1 | 1,0}$ \\
  \cline{2-3}
  & $\sfs{\rpp}{(2,1)\ (3,1)}$ &
    $H_{\twt_6^1 | 3,-2}$,
    $H_{\twt_6^1 | 3,-1}$,
    $H_{\twt_6^2 | 3,-2}$,
    $H_{\twt_6^2 | 3,-1}$,
    $H_{\twt_6^3 | 3,-1}$, \\
  & &
    $K_{\twt_5^1 | 3,-1}$,
    $K_{\twt_5^2 | 3,-2}$,
    $K_{\twt_5^2 | 3,-1}$,
    $K_{\twt_5^3 | 3,-2}$,
    $K_{\twt_5^3 | 3,-1}$ \\
  \cline{2-3}
  & $\sfs{\discref}{(2,1)\ (3,1)}$ &
    $H_{\twt_6^4 | 3,-1}$,
    $K_{\twt_5^4 | 3,-2}$,
    $K_{\twt_5^4 | 3,-1}$ \\
  \hline \hline

8 & $\torus \times I / \homtwo{3}{1}{1}{0}$ &
    $B_{T_6^2 | -3,1 | 1,0}$,
    $B_{T_6^2 | -2,3 | 1,-1}$,
    $B_{T_6^2 | -1,3 | 1,-2}$, \\
  & &
    $B_{T_7 | -2,-1 | -1,0}$,
    $B_{T_7 | -1,-1 | -2,-1}$,
    $B_{T_7 | 2,1 | 1,0}$, \\
  & &
    $B_{T_8^1 | -1,-1 | -1,0}$,
    $B_{T_8^1 | 1,1 | 1,0}$,
    $B_{T_8^2 | 0,1 | 1,1}$,
    $B_{T_8^2 | 1,1 | 1,0}$ \\
  \cline{2-3}
  & $\torus \times I / \homtwo{3}{2}{2}{1}$ &
    $B_{T_6^2 | -1,2 | 2,-3}$,
    $B_{T_7 | -1,-2 | -1,-1}$ \\
  \cline{2-3}
  & $\sfs{\rpp}{(2,1)\ (4,1)}$ &
    $H_{\twt_6^1 | 4,-3}$,
    $H_{\twt_6^1 | 4,-1}$,
    $H_{\twt_6^2 | 4,-3}$,
    $H_{\twt_6^2 | 4,-1}$,
    $H_{\twt_6^3 | 4,-1}$, \\
  & &
    $K_{\twt_5^1 | 4,-1}$,
    $K_{\twt_5^2 | 4,-3}$,
    $K_{\twt_5^2 | 4,-1}$,
    $K_{\twt_5^3 | 4,-3}$,
    $K_{\twt_5^3 | 4,-1}$ \\
  \cline{2-3}
  & $\sfs{\rpp}{(2,1)\ (5,2)}$ &
    $H_{\twt_6^1 | 5,-3}$,
    $H_{\twt_6^1 | 5,-2}$,
    $H_{\twt_6^2 | 5,-3}$,
    $H_{\twt_6^2 | 5,-2}$,
    $H_{\twt_6^3 | 5,-2}$, \\
  & &
    $K_{\twt_5^1 | 5,-2}$,
    $K_{\twt_5^2 | 5,-3}$,
    $K_{\twt_5^2 | 5,-2}$,
    $K_{\twt_5^3 | 5,-3}$,
    $K_{\twt_5^3 | 5,-2}$ \\
  \cline{2-3}
  & $\sfs{\rpp}{(3,1)\ (3,1)}$ &
    $H_{\twt_6^1 | 3,-1 | 3,-2}$,
    $H_{\twt_6^2 | 3,-1 | 3,-2}$,
    $H_{\twt_6^3 | 3,-1 | 3,-2}$, \\
  & &
    $K_{\twt_5^1 | 3,-1 | 3,-2}$,
    $K_{\twt_5^2 | 3,-2 | 3,-1}$,
    $K_{\twt_5^2 | 3,-1 | 3,-2}$,
    $K_{\twt_5^3 | 3,-1 | 3,-2}$ \\
  \cline{2-3}
  & $\sfs{\rpp}{(3,1)\ (3,2)}$ &
    $H_{\twt_6^1 | 3,-2 | 3,-2}$,
    $H_{\twt_6^1 | 3,-1 | 3,-1}$,
    $H_{\twt_6^2 | 3,-2 | 3,-2}$, \\
  & &
    $H_{\twt_6^2 | 3,-1 | 3,-1}$,
    $H_{\twt_6^3 | 3,-1 | 3,-1}$,
    $K_{\twt_5^1 | 3,-1 | 3,-1}$, \\
  & &
    $K_{\twt_5^2 | 3,-1 | 3,-1}$,
    $K_{\twt_5^3 | 3,-2 | 3,-2}$,
    $K_{\twt_5^3 | 3,-1 | 3,-1}$ \\
  \cline{2-3}
  & $\sfs{\discref}{(2,1)\ (4,1)}$ &
    $H_{\twt_6^4 | 4,-1}$,
    $K_{\twt_5^4 | 4,-3}$,
    $K_{\twt_5^4 | 4,-1}$ \\
  \cline{2-3}
  & $\sfs{\discref}{(2,1)\ (5,2)}$ &
    $H_{\twt_6^4 | 5,-2}$,
    $K_{\twt_5^4 | 5,-3}$,
    $K_{\twt_5^4 | 5,-2}$ \\
  \cline{2-3}
  & $\sfs{\discref}{(3,1)\ (3,1)}$ &
    $H_{\twt_6^4 | 3,-1 | 3,-1}$,
    $K_{\twt_5^4 | 3,-2 | 3,-1}$,
    $K_{\twt_5^4 | 3,-1 | 3,-2}$ \\
  \cline{2-3}
  & $\sfs{\discref}{(3,1)\ (3,2)}$ &
    $H_{\twt_6^4 | 3,-1 | 3,-2}$,
    $K_{\twt_5^4 | 3,-1 | 3,-1}$ \\
\hline
\end{tabular} \end{center}
\label{tab-trinames}
\end{table}
\renewcommand{\arraystretch}{1}

\vspace{1cm}
\noindent
Benjamin A.~Burton \\
Department of Mathematics, SMGS, RMIT University \\
GPO Box 2476V, Melbourne, 3010 VIC, Australia \\
(bab@debian.org)

\end{document}